\documentclass{siamltex}
\usepackage{amssymb}
\usepackage{graphicx}

\renewcommand{\thefootnote}{\fnsymbol{footnote}}

\begin{document}

\title{Pt\'ak's nondiscrete induction and its application to matrix iterations}

\author{J\"org Liesen\footnotemark[1]}

\footnotetext[1]{Institute of Mathematics, Technical University of Berlin,
Stra{\ss}e des 17. Juni 136, 10623 Berlin, Germany ({\tt liesen@math.tu-berlin.de}).}

\maketitle

\renewcommand{\thefootnote}{\arabic{footnote}}

\begin{abstract}
Vlastimil Pt\'ak's method of nondiscrete induction is based on the idea that
in the analysis of iterative processes one should aim at rates of convergence
as functions rather than just numbers, because functions may give convergence
estimates that are tight throughout the iteration rather than just asymptotically. 
In this paper we motivate and prove a theorem on nondiscrete induction
originally due to Potra and Pt\'ak,
and we apply it to the Newton iterations for computing the matrix polar decomposition
and the matrix square root. Our goal is to illustrate the application of the
method of nondiscrete induction in the finite dimensional numerical linear algebra
context. We show the sharpness of the resulting convergence estimate analytically
for the polar decomposition iteration and for special cases of the square
root iteration, as well as on some numerical examples for the square root iteration.
We also discuss some of the method's limitations and possible extensions.
\end{abstract}



\section{Introduction}
In the late 1960s, Vlastimil Pt\'ak (1925--1999) derived the {\em method of
nondiscrete induction}. This method for estimating the convergence of
iterative processes was originally motivated by a quantitative refinement
Pt\'ak had found for the closed graph theorem from functional
analysis in 1966~\cite{Pta66}. He published about 15 papers on the method,
five of them jointly with Potra in the 1980s, and the work on the method
culminated with their 1984 joint monograph~\cite{PotPtaBook84}. For historical
remarks on the development of the method see~\cite[Preface]{PotPtaBook84}
or~\cite[pp.~67--68]{Pta89}. 
Pt\'ak described the general motivation for the method of nondiscrete induction
in his paper ``What should be a rate of convergence?'', published in
1977~\cite{Pta77}:

\medskip
\begin{quote}
It seems therefore reasonable to look for another method of estimating the
convergence of iterative processes, one which would satisfy the following
requirements.
\begin{itemize}
\item[$1^\circ$] It should relate quantities which may be measured or estimated
during the actual process.
\item[$2^\circ$] It should describe accurately in particular the initial stage
of the process, not only its asymptotic behaviour since, after all, we
are interested in keeping the number of steps necessary to obtain a good
estimate as low as possible.
\end{itemize}
\end{quote}

\medskip
This seems to be almost too much to ask for. Yet, for some iterations the
method of nondiscrete induction indeed leads to analytical convergence
estimates which satisfy the above requirements. The essential idea, as we will
describe in more detail below, is that in the method of nondiscrete induction
the rate of convergence of an iterative process is considered {\em a function rather
than just a number}. Pt\'ak derived convergence results satisfying the above
requirements for a handful of examples including Newton's method~\cite{Pta76a}
(this work was refined in~\cite{PotPta80}), and an iteration for solving a certain
eigenvalue problem~\cite{Pta76b}.
For these examples it was shown that the convergence estimates
resulting from the method of nondiscrete induction are indeed
optimal in the sense that in certain cases they are attained
in every step of the iteration. In addition to
the original papers, comprehensive statements of such sharpness results are
given in the book of Potra and Pt\'ak; see,
e.g.,~\cite[Proposition~5.10]{PotPtaBook84} for Newton's method
and~\cite[Proposition~7.5]{PotPtaBook84} for the eigenvalue iteration.

Despite these strong theoretical results, it appears that the method of nondiscrete induction
never became widely known. Even in the literature on Newton methods for nonlinear
problems it is often mentioned only marginally (if at all); see, e.g., Deuflhard's
monograph~\cite[p.~49]{DeuBook11}. Part of the reason for this neglect of the
method may be the lack of numerically computed examples in the original publications.

The goals of this paper are: (1) to motivate and prove a theorem on nondiscrete
induction due to Potra and Pt\'ak that is directly applicable in the analysis of
iterative processes, (2) to explain on two examples (namely the Newton iterations 
for the matrix polar decomposition and the matrix square root) how the method 
of nondiscrete induction can be applied to matrix iterations in numerical linear algebra, 
and (3) to demonstrate the method's effectiveness as well as discuss some of its 
weaknesses in the numerical linear algebra context.
It must be stressed upfront, that most theoretical results in Sections~\ref{sec:polar}
and~\ref{sec:appl} also could be derived using the general theory of Newton's
method for nonlinear operators in Banach spaces as
described in~\cite[Chapters~2 and~5]{PotPtaBook84} and in the related papers of
Potra and Pt\'ak, in particular~\cite{PotPta80,Pta76a}. The strategy in this paper is,
however, to apply the method of nondiscrete induction without any functional analytic
framework and differentiability assumptions directly to the given algebraic matrix
iterations.

The paper is organized as follows. In Section~\ref{sec:method} we describe the method
of nondiscrete induction and derive the theorem of Potra and Pt\'ak.
We then apply this theorem in analysis the Newton iterations for computing
the matrix polar decomposition (Section~\ref{sec:polar})
and the matrix square root (Section~\ref{sec:appl}). For both iterations we prove
convergence results and illustrate them numerically. In Section~\ref{sec:concl} we give
concluding remarks and an outlook to further work.

\section{The method of nondiscrete induction}\label{sec:method}

In most of the publications on the method of nondiscrete induction, Pt\'ak
formulated the method based on his ``Induction Theorem'', which is an inclusion
result for certain subsets of a metric space; see,
e.g.,~\cite[p.~280]{Pta76a}, \cite[p.~225]{Pta76b}, \cite[p.~282]{Pta77},
\cite[p.~52]{Pta89},
or the book of Potra and Pt\'ak~\cite[Proposition~1.7]{PotPtaBook84}.
Instead of the original Induction Theorem we will in this paper use a
result that was stated and proven by Potra and Pt\'ak as a ``particular
case of the Induction Theorem'' in~\cite[Proposition~1.9]{PotPtaBook84};
also cf.~\cite[p.~66]{PotPta80}. Unlike the Induction Theorem,
this result is formulated directly in terms of an iterative algorithm, and
therefore it is easier to apply in our context. In order to fully explain
the result and its consequences for the analysis of iterative algorithms
we will below give a motivation of the required concepts as well as a
complete proof of the assertion without referring to the Induction Theorem.

Let $(E,d)$ be a complete metric space, where $d(x,y)$ denotes the distance
between $x,y\in E$. Consider a mapping $G\,:\,D(G)\rightarrow E$,
where $D(G)\subseteq E$ denotes the domain of definition of $G$, i.e., the subset of all
$x\in E$ for which $G(x)\in E$ is well defined. For each $x_0\in D(G)$ we
may then consider the iterative algorithm $(G,x_0)$ with iterates given by
\begin{equation}\label{eqn:algorithm}
x_{k+1} = G(x_k),\quad k=0,1,2,\dots\,.
\end{equation}
If all iterates are well defined, i.e., $G(x_k)\in D(G)$ for all $k\geq 0$, then
the iterative algorithm $(G,x_0)$ is called {\em meaningful}.

It is clear from (\ref{eqn:algorithm}) that a (meaningful) iterative algorithm can converge to
some $x_*\in E$ only when the distances between $x_k$ and $x_{k+1}=G(x_k)$ decrease to
zero for $k\rightarrow \infty$. This requirement, and the rate by which the convergence to zero
occurs, are formalized in the method of nondiscrete induction by means of a (nondiscrete)
family of subsets $Z(t)\subseteq D(G)$ that depend on a positive real parameter $t$. The
key idea is that for each $t\in (0,\tau)$, where possibly $\tau=\infty$, the set $Z(t)$
should contain all elements $x\in D(G)$ for which the distance between $x$ and $G(x)$ is at
most $t$, i.e. $d(x,G(x))\leq t$, and then to analyze what happens with the elements of
$Z(t)$ under one application of the mapping $G$. This means that the method compares the
values
$$d(x,G(x))\quad\mbox{and}\quad d(G(x),G(G(x))).$$
Hopefully, the new distance $d(G(x),G(G(x)))$ is (much) smaller than the previous distance
$d(x,G(x))$. Using the family of sets this can formally be written as
$$G(x)\in Z(\omega(t)),$$
where $\omega$ should be a (positive) ``small function'' of $t$ that needs to be determined.
This function is called the {\em rate of convergence} of the iterative algorithm. Since in
the limit the distances between the iterates should approach zero, the limiting point(s)
$x_*$ of the algorithm are contained in the limit set $Z(0)$, which is not to be constructed
explicitly but is defined as
\begin{equation}\label{eqn:Z0}
Z(0)\equiv \lim_{t\rightarrow 0} Z(t)= \bigcap_{0<t<\tau}\,\bigl(\overline{\bigcup_{s\leq t} Z(s)}\bigr)
\;\subseteq\; \overline{D(G)},
\end{equation}
where $\overline{S}$ denotes the closure of the set $S$. Note that the set $Z(0)$ possibly is not
a subset of $D(G)$.

We need the following formal definition of a rate of convergence.

\medskip
\begin{definition}\label{def:rate}
Let $T=(0,\tau)$ be a given real interval, where possibly $\tau=\infty$. Let $\omega\,:\,T\rightarrow T$
be a function and denote by $\omega^{(k)}$ its $k$th iterate, i.e., $\omega^{(0)}(t)=t$ and
$\omega^{(k)}(t)=\omega(\omega^{(k-1)}(t))$ for $k\geq 1$. The function $\omega$ is called 
a rate of convergence on $T$ when the corresponding series
\begin{equation}\label{eqn:sigma}
\sigma(t)\equiv t+\omega(t)+\omega^{(2)}(t)+\dots = \sum_{k=0}^\infty \omega^{(k)}(t)
\end{equation}
converges for all $t\in T$.
\end{definition}

\medskip
The reason for the convergence of the series $\sigma$ in~(\ref{eqn:sigma}) will become apparent
in the proof of the following result of Potra and Pt\'ak, cf.~\cite[Proposition~1.9]{PotPtaBook84}.

\medskip
\begin{theorem}\label{thm:PotPta}
Let $(E,d)$ be a complete metric space, let $G\,:\,D(G)\rightarrow E$ be a given mapping, and
let $x_0\in D(G)$. If there exist a rate of convergence $\omega$ on a real interval $T=(0,\tau)$,
with the corresponding function $\sigma$ as in (\ref{eqn:sigma}), and a family of sets
$Z(t)\subseteq D(G)$ for all $t\in T$, such that
\begin{itemize}
\item[(i)] $x_0\in Z(t_0)$ for some $t_0\in T$,
\item[(ii)] $d(x,G(x))\leq t$ and $G(x)\in Z(\omega(t))$ for each $t\in T$ and $x\in Z(t)$,
\end{itemize}
then:
\begin{enumerate}
\item[(1)] $(G,x_0)$ is a meaningful iterative algorithm and the sequence $\{x_k\}$ of its iterates
converges to a point $x_*\in Z(0)\subseteq \overline{D(G)}$. (In particular, $Z(0)\neq \emptyset$.)
\item[(2)] The following relations are satisfied for all $k\geq 0$:\\
a) $x_k\in Z(\omega^{(k)}(t_0))$,\\
b) $d(x_k,x_{k+1})\leq \omega^{(k)}(t_0)$,\\
c) $d(x_k,x_*)\leq \sigma(\omega^{(k)}(t_0))$.\\
The right hand side of the bound in (2.c) is a strictly monotonically decreasing
function of $k$. Moreover, if equality holds in (2.c) for some $k_0\geq 0$, then
equality holds for all $k\geq k_0$.
\end{enumerate}
\end{theorem}

\medskip
\begin{proof}
Since $x_0\in Z(t_0)\subseteq D(G)$, we know that $x_1=G(x_0)$ exists. Now (ii) implies that
$$x_1\in Z(\omega(t_0))\quad\mbox{and}\quad d(x_0,x_1)\leq t_0.$$
We can apply the same reasoning to $x_1\in Z(\omega(t_0))\subseteq D(G)$, which yields the
existence of $x_2=G(x_1)$ with
$$x_2\in Z(\omega^{(2)}(t_0))\quad\mbox{and}\quad d(x_1,x_2)\leq \omega(t_0).$$
Continuing in this way we obtain a sequence of well defined iterates $x_{k+1}=G(x_{k})\in D(G)$ with
$$x_k\in Z(\omega^{(k)}(t_0))\quad\mbox{and}\quad
d(x_{k},x_{k+1})\leq \omega^{(k)}(t_0),\quad k\geq 0.$$
This shows items (2.a) and (2.b).

Next observe that, for all $k\geq 0$ and $m\geq 1$,
\begin{eqnarray*}
d(x_{k},x_{k+m})&\leq& d(x_{k},x_{k+1})+\dots+d(x_{k+m-1},x_{k+m})\\
&\leq& \omega^{(k)}(t_0)+\dots+\omega^{(k+m-1)}(t_0),
\end{eqnarray*}
where we have used (2.b). Convergence of the series $\sigma(t)$ for all $t\in T$
implies that the sequence $\{x_k\}\subset D(G)$ is a Cauchy sequence, which by
completeness of the space converges to a well defined limit point $x_*$. From
$x_k\in Z(\omega^{(k)}(t_0))$ and $\omega^{(k)}(t_0)\rightarrow 0$ for
$k\rightarrow\infty$, we obtain $x_*\in Z(0)\subseteq \overline{D(G)}$. Hence we
have shown (1).

For all $k\geq 0$ and $m\geq 1$,
\begin{eqnarray*}
d(x_{k},x_{*})&\leq& d(x_{k},x_{k+1})+\dots+d(x_{k+m-1},x_{k+m})+d(x_{k+m},x_*)\\
&\leq& \omega^{(k)}(t_0)+\dots+\omega^{(k+m-1)}(t_0)+d(x_{k+m},x_*)\\
&=&d(x_{k+m},x_*) +\sigma_{m-1}(\omega^{(k)}(t_0)),\quad \mbox{where}\;\;
\sigma_{m-1}(t)\equiv\sum_{j=0}^{m-1}\omega^{(j)}(t).
\end{eqnarray*}
Taking the limit $m\rightarrow\infty$ shows that $d(x_{k},x_{*})\leq
\sigma(\omega^{(k)})(t_0)$, which proves (2.c).

Since $\omega:T\rightarrow T$ and $t_0\in T$, we have $\omega^{(k)}(t_0)>0$ for all $k\geq 0$,
and thus
$$\sigma(\omega^{(k+1)}(t_0))=\sum_{j=k+1}^\infty \omega^{(j)}(t_0)=
\left(\sum_{j=k}^\infty \omega^{(j)}(t_0)\right)-\omega^{(k)}(t_0)<\sigma(\omega^{(k)}(t_0)).$$
For the proof of the last assertion about the equality in (2.c) we refer to the proof
of~\cite[Proposition~1.11]{PotPtaBook84}. (In that proof it should read ``If (10) is
verified'' instead of ``If (9) is verified''.)
\end{proof}

\medskip
From the formulation of Theorem~\ref{thm:PotPta} it becomes apparent why Pt\'ak
considered his method a ``continuous analogue'' of the classical mathematical induction;
also cf. his own explanation in~\cite[p.~225--226]{Pta76b}: In condition (i) we require
the existence of at least one set $Z(t_0)$ that contains the initial value $x_0$.
This corresponds to the base step in the classical mathematical induction.
Condition (ii) then considers what happens with the elements of the set $Z(t)$
under one application of the mapping $G$. This corresponds to the inductive step.

In the classical theory of iterative processes, an iteration $x_{k+1}=G(x_k)$
is called {\em convergent of order $p\geq 1$} when there exist a constant $c>0$ and a
positive integer $k_0$, such that
\begin{equation}\label{eqn:classical}
d(x_*,x_{k+1})\,\leq\, c\, d(x_*,x_{k})^p\quad\mbox{for all $k\geq k_0$.}
\end{equation}
In particular, $p=1$ and $p=2$ give ``linear'' and ''quadratic'' convergence,
respectively, where in the first case we require that $c\in (0,1)$.
As Pt\'ak pointed out, e.g., in~\cite{Pta77,Pta89}, such classical bounds compare quantities that
are not available at any finite step of the process, simply because $x_*$ is unknown. In addition,
since $p$ is the same fixed number for all $k\geq k_0$ and $k_0$ often is large, such
bounds in many cases cannot give a good description of the initial stages of the iteration;
cf. the quote from~\cite{Pta77} given in the Introduction above.
Item (2.c) in Theorem~\ref{thm:PotPta}, on the other hand, gives an {\em a priori} bound
on the error norm in each step $k\geq 0$ of the iterative algorithm. Moreover, since the
convergence rate is a function (rather than just a number), there appears to be a better chance
that the bound is tight throughout the iteration.

Finally, we note that from the statement of Theorem~\ref{thm:PotPta} it is clear
that if $x_*\in D(G)$ and $G$ is continuous at $x_*$, then $x_*=G(x_*)$,
i.e., $x_*$ is a fixed point of $G$. In contrast to classical fixed point
results, however, Theorem~\ref{thm:PotPta} contains no explicit continuity assumption on~$G$.
As shown in~\cite[pp.~10--12]{PotPtaBook84} (see also~\cite[Section~3]{Pta76b}), 
the Banach fixed point theorem, which assumes Lipschitz continuity of $G$, can be easily derived
as a corollary of Theorem~\ref{thm:PotPta}.

\section{Application to the Newton iteration for the matrix polar decomposition}\label{sec:polar}

Let $A\in {\mathbb C}^{n\times n}$ be nonsingular, and let $A=W\Sigma V^*$ be a singular value
decomposition with $\Sigma={\rm diag}(\sigma_1,\dots,\sigma_n)$, $\sigma_j>0$ for $j=1,\dots,n$,
and with unitary matrices $W,V\in {\mathbb C}^{n\times n}$. The factorization
$$A=UH,\quad\mbox{where}\quad U\equiv WV^*,\quad H\equiv V\Sigma V^*$$
is called a polar decomposition of $A$. Here $U$ is unitary and $H$ is Hermitian positive definite.
The theory of the matrix polar decomposition and many algorithms for computing this decomposition are
described in Higham's monograph on matrix functions~\cite[Chapter~8]{HigBook08}.

One of the algorithms for computing the matrix polar decomposition is Newton's
iteration~\cite[Equation (8.17)]{HigBook08}:
\begin{eqnarray}
& & X_0=A.\label{eqn:polar1}\\
& & X_{k+1} = \frac12 (X_k+X_k^{-*}),\quad k=0,1,2,\dots\,.\label{eqn:polar2}
\end{eqnarray}
This iteration can be derived by applying Newton's method to the matrix equation
$X^*X-I=0$; cf.~\cite[Problem~8.18]{HigBook08}. As shown in~\cite[Theorem~8.12]{HigBook08},
its iterates $X_k$ converge to the unitary polar factor $U$, and they satisfy the
convergence bound
$$\|U-X_{k+1}\|\leq \frac12\,\|X_k^{-1}\|\,\|U-X_k\|^2,\quad k=0,1,2,\dots\,,$$
where $\|\cdot\|$ denotes any unitarily invariant and submultiplicative matrix norm.
This is a bound of the form (\ref{eqn:classical}) that shows quadratic convergence of
the iteration (\ref{eqn:polar1})--(\ref{eqn:polar2}). The derivation of this bound
as well as further a posteriori convergence bounds for this iteration in terms of the
singular values of $A$ can also be found in~\cite[Theorem~3.1]{Hig86b}.

\subsection{A convergence theorem for the iteration (\ref{eqn:polar1})--(\ref{eqn:polar2})}
We will now analyze the algebraic matrix iteration (\ref{eqn:polar1})--(\ref{eqn:polar2})
using Theorem~\ref{thm:PotPta}. In the notation of the theorem we consider
$E={\mathbb C}^{n\times n}$ and we write the distance between $X,Y\in E$ as $\|X-Y\|$,
where $\|\cdot\|$ denotes the spectral norm. We also consider the mapping $G\,:\,D(G)\rightarrow E$
defined by
$$G(X)\equiv\frac12 (X+X^{-*});\quad\mbox{cf. (\ref{eqn:polar2}).}$$
We need to determine $D(G)$ and define the sets $Z(t)\subseteq D(G)$ for some
interval $T=(0,\tau)$ on which the nondiscrete induction is based. The value of $\tau$
is not specified yet; it will be determined when the rate of convergence
has been found.

Obviously, $D(G)$ is equal to the set of the invertible matrices in
${\mathbb C}^{n\times n}$, and hence every $X\in Z(t)$ must be invertible.
We also require that every $X\in Z(t)$ satisfies $\|X-G(X)\|\leq t$; cf.
the motivation of Theorem~\ref{thm:PotPta} and the first assumption in (ii)
of the theorem. Moreover, for $X_0=A=W\Sigma V^*$ we get
$X_1=G(X_0)=\frac12 W(\Sigma+\Sigma^{-1})V^*$.
Inductively we see that every $X_k$ that can possibly occur in the iteration
(\ref{eqn:polar1})--(\ref{eqn:polar2}) is of the form $X_k=WSV^*$ for some diagonal
matrix $S$ with real and positive diagonal elements. For the given matrix $A$ and
each $t\in T$ we therefore define
\begin{eqnarray*}
Z(t) \equiv \bigl\{  X\in {\mathbb C}^{n\times n} & | &\; X=WSV^*,\;\;S={\diag}(s_1,\dots,s_n),\;\;
s_j>0,\;\; j=1,\dots,n, \\
& & \;\|X-G(X)\|\leq t \bigr\}.
\end{eqnarray*}
In particular, $Z(t)\subset D(G)$.

Now let $X\in Z(t)$, so that $X=WSV^*$ for some $S={\diag}(s_1,\dots,s_n)$ with $s_j>0$ for $j=1,\dots,n$.
Then $\widehat{X}\equiv G(X)=\frac12 W(S+S^{-1})V^*$, and
\begin{eqnarray}
& & X-\widehat{X} = \frac12 W(S-S^{-1})V^*,\nonumber\\
& & \|X-\widehat{X}\| = \max_{1\leq j\leq n} f(s_j),\quad\mbox{where}\quad
f(s)\equiv \frac12\,|s-s^{-1}|.\label{eqn:pp1}
\end{eqnarray}
From (\ref{eqn:pp1}) we see that in the limit $t\rightarrow 0$ we
obtain $Z(0)=\{WV^*\}$, i.e., the limiting set $Z(0)$ consists only of the unitary polar factor of $A$.

In order to determine a rate of convergence we compute
\begin{eqnarray}
\widehat{X}-G(\widehat{X}) &=& \frac12 \widehat{X}^{-*}(\widehat{X}^*\widehat{X}-I)
= W (S+S^{-1})^{-1}\,\frac14\,(S-S^{-1})^2W^*,\label{eqn:eq}
\end{eqnarray}
and hence
\begin{eqnarray}
\|\widehat{X}-G(\widehat{X})\| &=& \max_{1\leq j\leq n}g(s_j),\quad\mbox{where}\label{eqn:pp2}\\
g(s) &\equiv& \frac{\frac14 (s-s^{-1})^2}{s+s^{-1}}=\frac{f^2(s)}{2(f^2(s)+1)^{1/2}}.\label{eqn:pp3}
\end{eqnarray}
Here we have used the definition of $f(s)$ from (\ref{eqn:pp1}) and that
$(s+s^{-1})^2=4(f^2(s)+1)$. Note that $g(1)=f(1)=0$ and $g(s)<f(s)$ for
all positive $s\neq 1$.

\medskip
\begin{lemma}
The maxima in (\ref{eqn:pp1}) and (\ref{eqn:pp2}) are attained at the same point $s_i$.
\end{lemma}

\medskip
\begin{proof}
Without loss of generality we can assume that $0<s_1\leq s_2\leq\cdots\leq s_n$.
Elementary computations show that
$$\mbox{$f'(s)<0$ for $0<s<1\;\;$ and $\;\;f'(s)>0$ for $s>1$,}$$
$$\mbox{$g'(s)<0$ for $0<s<1\;\;$ and $\;\;g'(s)>0$ for $s>1$,}$$
i.e., both $f$ and $g$ are strictly monotonically decreasing for $0<s<1$ and
strictly monotonically increasing for $s>1$. This means that
$$\max_{1\leq j\leq n} f(s_j)=\max\{f(s_1),f(s_n)\}\quad \mbox{and}\quad
\max_{1\leq j\leq n} g(s_j)=\max\{g(s_1),g(s_n)\}.$$

If $0<s_1\leq s_n\leq 1$ resp. $1\leq s_1\leq s_n$, then obviously the maximum for both
functions is attained at $s_1$ resp. $s_n$.
In the remaining case $0<s_1<1<s_n$ we will use
that $f(s)=f(s^{-1})$ and $g(s)=g(s^{-1})$ for all $s>0$, which can be easily shown.
Suppose that the maximum in (\ref{eqn:pp1}) is attained at $s_n$, i.e., $f(s_1)\leq f(s_n)$.
Then also $f(s_1^{-1})\leq f(s_n)$, and thus we must have $s_1^{-1}\leq s_n$ because of the
(strict) monotonicity of $f$ for $s>1$. But since $g$ is also (strictly) monotonically
increasing for $s>1$, we get $g(s_1)=g(s_1^{-1})\leq g(s_n)$. An analogous argument
applies when the maximum in (\ref{eqn:pp1}) is attained at $s_1$.
\end{proof}

\medskip
For $X=WSV^*\in Z(t)$ there exist an index $i$
and a real number $\widetilde{t}>0$ such that
$$\|X-\widehat{X}\| = \max_{1\leq j\leq n} f(s_j) = f(s_i) = \widetilde{t} \leq t.$$
The previous lemma then implies
\begin{eqnarray}
\|\widehat{X}-G(\widehat{X})\| &=& \max_{1\leq j\leq n} g(s_j) = g(s_i) =
\frac{f^2(s_i)}{2(f^2(s_i)+1)^{1/2}}=
\frac{\widetilde{t}^2}{2(\widetilde{t}^2+1)^{1/2}}\label{eqn:exact2}\\
&\leq &\frac{t^2}{2(t^2+1)^{1/2}}\equiv \omega(t),\label{eqn:omegaPP}
\end{eqnarray}
where we have used that the function $\omega$ introduced in (\ref{eqn:omegaPP}) is (strictly)
monotonically increasing for $t\in T\equiv (0,\infty)$.
Pt\'ak showed in~\cite{Pta76a} (also cf.~\cite{PotPta80,Pta77}
or~\cite[Example~1.4]{PotPtaBook84})
that $\omega$ is a rate of convergence on $T$ in the sense of Definition~\ref{def:rate}
with the corresponding series $\sigma$ given by\footnote{Here the notation $\sigma$ should
not be confused with the singular values $\sigma_j$ of $A$. This slight inconvenience
is the price for retaining the original notation of Pt\'ak~\cite{Pta76a,Pta76b,Pta77}.}
\begin{equation}\label{eqn:sigmaPP}
\sigma(t)=t-1+(t^2+1)^{1/2}.
\end{equation}
He pointed out in~\cite{Pta77} that $\omega(t)\approx t/2$ for large $t$, and $\omega(t)\approx t^2/2$
for small $t$. Thus, the rate of convergence $\omega$ can describe both a linear and a quadratic
phase of the convergence; cf. the numerical examples in Section~\ref{sec:numerPP}.

We have shown above that $X\in Z(t)$ implies $G(X)\in Z(\omega(t))$ for $\omega$ defined in
(\ref{eqn:omegaPP}), and hence we have satisfied condition (ii) in Theorem~\ref{thm:PotPta}.
In order to satisfy condition (i) we simply set $t_0\equiv \|A-G(A)\|$, then $A\in Z(t_0)$
is guaranteed (recall that $X_0=A$). Theorem~\ref{thm:PotPta} now gives the following result.

\medskip
\begin{theorem}\label{thm:pp}
Let $A=W\Sigma V^*\in {\mathbb C}^{n\times n}$ with $\Sigma={\diag}(\sigma_1,\dots,\sigma_n)$,
$\sigma_j>0$ for  $j=1,\dots,n$, and unitary matrices $W,V\in {\mathbb C}^{n\times n}$.
Then the iterative algorithm (\ref{eqn:polar1})--(\ref{eqn:polar2}) is meaningful and its
sequence of iterates $\{X_k\}$ converges to the unitary polar factor $U=WV^*$ of $A$. Moreover,
with $t_0\equiv \frac12\max_{1\leq j\leq n} |\sigma_j-\sigma_j^{-1}|$,
and $\omega$ resp. $\sigma$ defined as in (\ref{eqn:omegaPP}) resp. (\ref{eqn:sigmaPP}), we have
\begin{equation}\label{eqn:boundPP}
\|U-X_k\|\leq \sigma(\omega^{(k)}(t_0)),\quad k=0,1,2,\dots\,.
\end{equation}
The inequality (\ref{eqn:boundPP}) is an equality for all $k\geq 1$.
\end{theorem}

\medskip
\begin{proof}
It remains to show that (\ref{eqn:boundPP}) is an equality for all $k\geq 1$. According to the
last part of Theorem~\ref{thm:PotPta} it suffices to show that $\|U-X_1\|=\sigma(\omega(t_0))$.
With $X_0=A=W\Sigma V^*$ we get $X_1=\frac12 W(\Sigma+\Sigma^{-1})V^*$, so that
\begin{eqnarray*}
\|U-X_1\|&=&\left\|I-\frac12 (\Sigma+\Sigma^{-1})\right\|=\max_{1\leq j\leq n}\,
\left|1-\frac12 (\sigma_j+\sigma_j^{-1})\right|\\
&=& \max_{1\leq j\leq n}\,\frac12 (\sigma_j+\sigma_j^{-1})-1
=\max_{1\leq j\leq n}\,\left(\frac14 (\sigma_j-\sigma_j^{-1})^2+1\right)^{1/2}-1\\
&=& (t_0^2+1)^{1/2}-1 = \sigma(t_0)-t_0 = \sigma(\omega(t_0)).
\end{eqnarray*}
(Recall that $\sigma(t)-t=\sigma(\omega(t))$ for all $t\in T$ by definition of the series $\sigma$.)
\end{proof}

\medskip
The fact that (\ref{eqn:boundPP}) is an equality for all $k\geq 1$ clearly demonstrates the advantage 
of using functions rather than just numbers for describing the convergence of iterative algorithms. 
The equality could be expected when considering that Potra and Pt\'ak showed that the convergence 
estimates obtained by the method of nondiscrete induction applied to Newton's method for nonlinear 
operators in Banach spaces are sharp for special cases of scalar functions; see,
e.g.,~\cite[Proposition~5.10]{PotPtaBook84}. One can interpret the iteration
(\ref{eqn:polar1})--(\ref{eqn:polar2}) in terms of $n$ uncoupled scalar iterations, and can then
derive Theorem~\ref{thm:pp} using results in~\cite[Chapter~5]{PotPtaBook84} or the original
papers~\cite{PotPta80,Pta76b}.

\subsection{Numerical examples}\label{sec:numerPP}

\begin{figure}
\begin{minipage}{0.49\linewidth}
\includegraphics[width=0.99\linewidth]{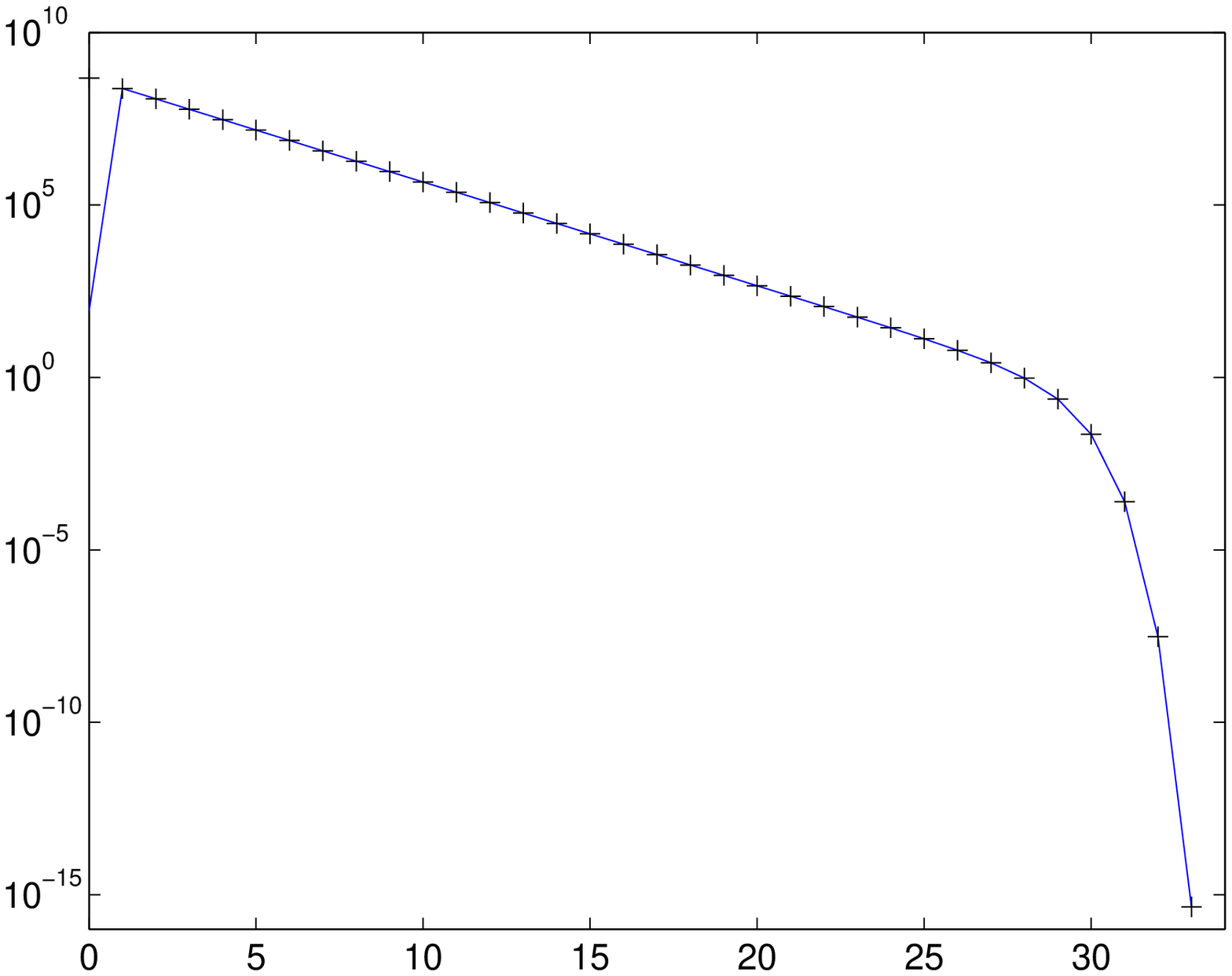}
\end{minipage}
\hfill
\begin{minipage}{0.49\linewidth}
\includegraphics[width=0.99\linewidth]{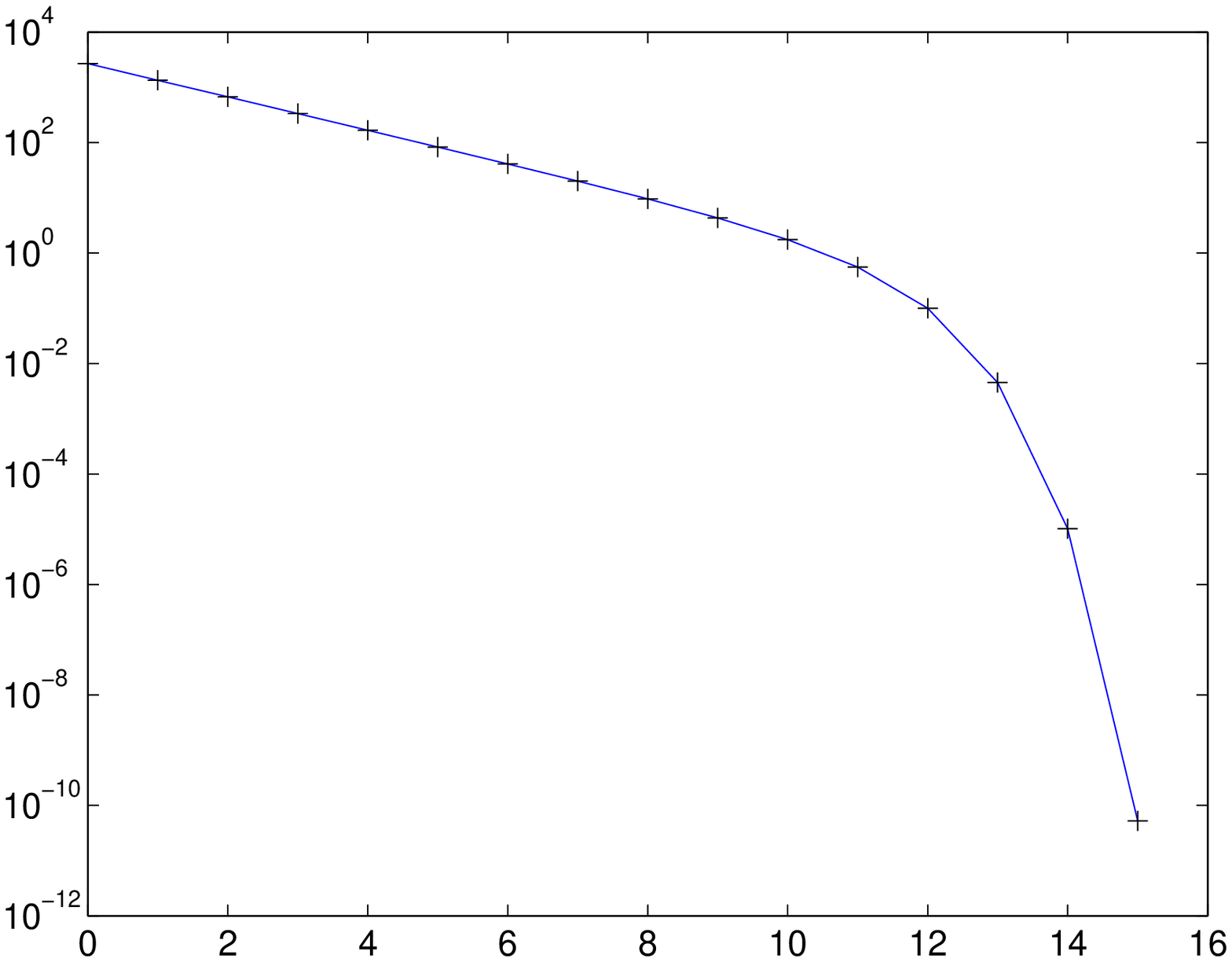}
\end{minipage}
\caption{$\|I-X_k\|$ (solid) and $\sigma(\omega^{(k)}(t_0))$ (pluses) for
the Moler matrix (left) and the Fielder matrix (right).}\label{fig:1pp}
\end{figure}

\begin{figure}
\begin{minipage}{0.49\linewidth}
\includegraphics[width=0.99\linewidth]{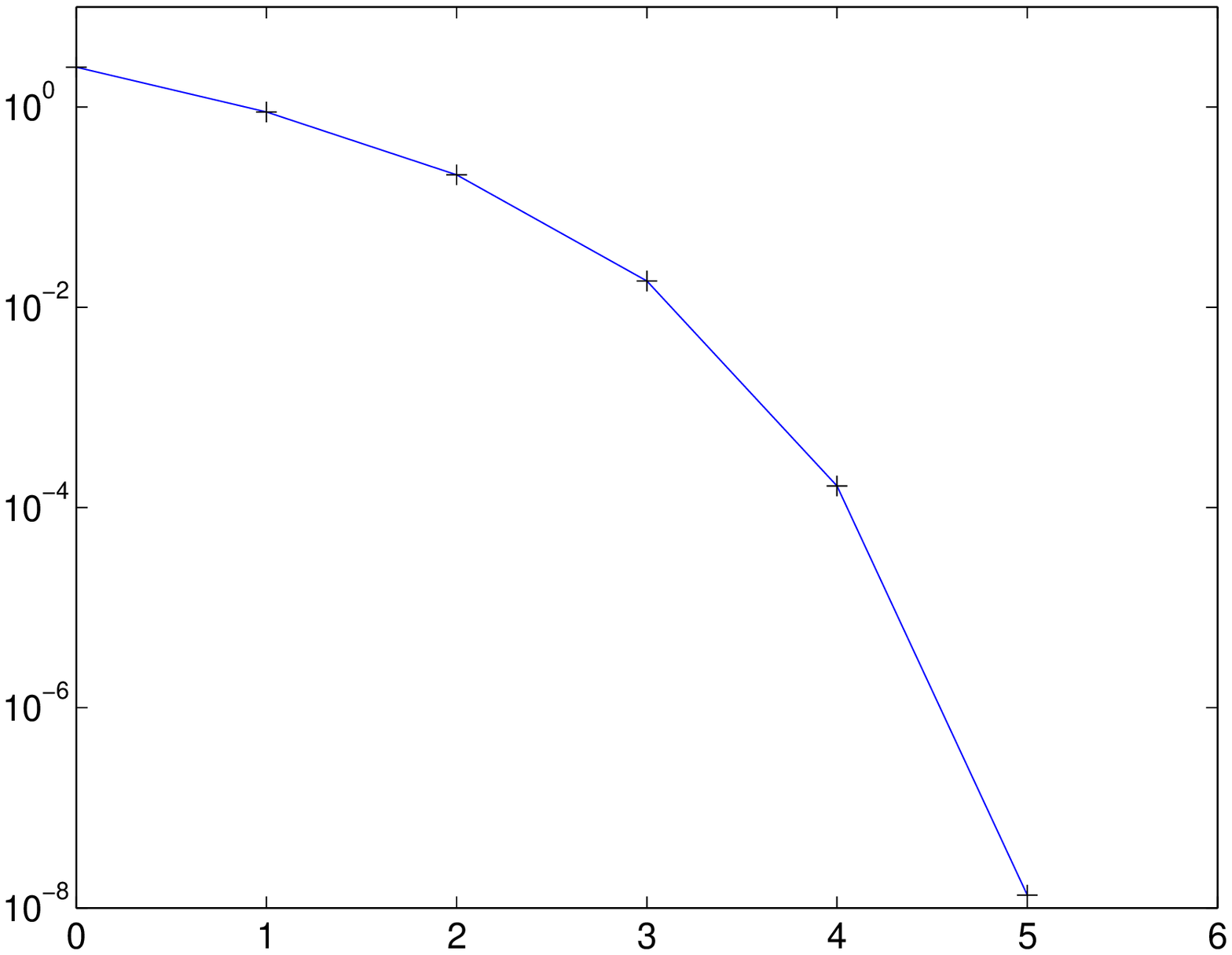}
\end{minipage}
\hfill
\begin{minipage}{0.49\linewidth}
\includegraphics[width=0.99\linewidth]{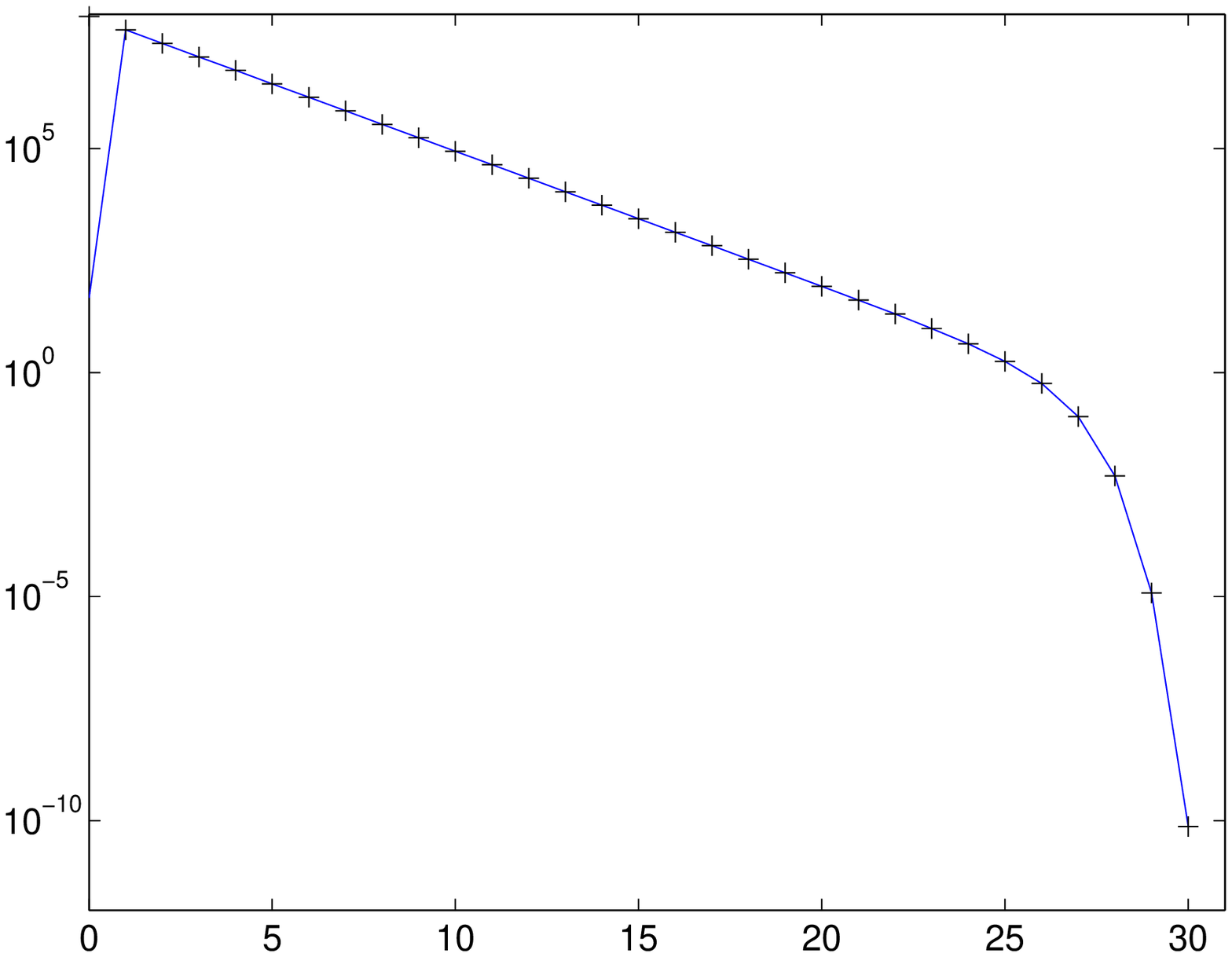}
\end{minipage}
\caption{$\|I-X_k\|$ (solid) and $\sigma(\omega^{(k)}(t_0))$ (pluses) for
the Jordan blocks matrix (left) and the Frank matrix (right).}\label{fig:2pp}
\end{figure}

We will illustrate the convergence of the iteration (\ref{eqn:polar1})--(\ref{eqn:polar2})
and the bound (\ref{eqn:boundPP}) on four examples. All experiments shown in this paper
were computed in MATLAB. For numerical purposes we first computed a singular value
decomposition $A=W\Sigma V^*$ using MATLAB's command {\tt svd} and then applied
(\ref{eqn:polar1})--(\ref{eqn:polar2}) with $X_0=\Sigma$. The error in every step $k\geq 0$,
given by $\|I-X_k\|$, is shown by the solid lines in Figures~\ref{fig:1pp} and~\ref{fig:2pp}.
The pluses ($+$) in these figures indicate the corresponding values of the
right hand side of (\ref{eqn:boundPP}). It should be stressed that the right hand side of
(\ref{eqn:boundPP}) is a simple scalar function that only depends on the value $t_0$.
This function may also be computed by the more explicit expressions given in~\cite{PotPta80}
or~\cite[Chapter~5]{PotPtaBook84}.

We present numerical results for four different matrices:

\smallskip
{\em Moler matrix.\/}
Generated using MATLAB's command {\tt gallery('moler',16)}, this is a $16\times 16$ symmetric
positive definite matrix with 15 eigenvalues between $2$ and $88$,
and one small eigenvalue of order $10^{-9}$.
For this matrix $t_0=2.3861e+08$; the numerical results are shown in the left
part of Figure~\ref{fig:1pp}.

\smallskip
{\em Fiedler matrix.\/}
Generated using {\tt gallery('fiedler',88)}, this is a $88\times 88$ symmetric indefinite matrix
with 87 negative eigenvalues and 1 positive eigenvalue.
For this matrix $t_0=1.3450e+03$; the numerical results are shown in
the right part of Figure~\ref{fig:1pp}.

\smallskip
{\em Jordan blocks matrix.\/} The $100\times 100$ matrix ${\rm diag}(J_{50}(1.5),J_{50}(2.5))$,
where $J_n(\lambda)$ is an $n\times n$ Jordan block with eigenvalue $\lambda$.
For this matrix $t_0=1.6064$; the numerical results are shown in the left part of Figure~\ref{fig:2pp}.

\smallskip
{\em Frank matrix.\/} Generated using {\tt gallery('frank',12)}, this is a $12\times 12$
nonnormal matrix with ill-conditioned real eigenvalues; the condition number of its eigenvector
matrix computed in MATLAB is $1.22e+8$. For this matrix
$t_0=4.4698e+07$; the numerical results are shown in the right part of Figure~\ref{fig:2pp}.

\smallskip
In the examples we observe the typical behavior of Newton's method: An initial phase of
linear convergence being followed by a phase of quadratic convergence in which the error is reduced
very quickly. While the convergence bound (\ref{eqn:boundPP}) is an equality for all $k\geq 1$,
it may happen that $\|I-X_0\|\ll \sigma(t_0)$ as well as $\|I-X_0\|\ll \|I-X_1\|$;
see the results for the Moler matrix and the Frank matrix. As shown in Theorem~\ref{thm:PotPta},
the error is guaranteed to decrease strictly monotonically after the first step. In practice
the iteration can be scaled for stability as well as for accelerating the convergence;
see~\cite{ByeXu08,Hig86b} or~\cite[Section~8.6]{HigBook08} for details. We give some brief
comments on scaled iterations in Section~\ref{sec:concl}.

\section{Application to the Newton iteration for the matrix square root}\label{sec:appl}

Let $A\in {\mathbb C}^{n\times n}$ be a given, possibly singular,
matrix. Any matrix $X\in {\mathbb C}^{n\times n}$
that satisfies $X^2=A$ is called a square root of $A$. For the theory of matrix square roots and
many algorithms for computing them we refer to~\cite[Chapter~6]{HigBook08}.

One of the algorithms for computing matrix square roots is Newton's
iteration; cf.~\cite[Equation (6.12)]{HigBook08}:
\begin{eqnarray}
& & \mbox{$X_0$ is invertible and commutes with $A$.}\label{eqn:newton1}\\
& & X_{k+1} = \frac12 (X_k+X_k^{-1}A),\quad k=0,1,2,\dots\,.\label{eqn:newton2}
\end{eqnarray}
Similar to the Newton iteration for the matrix polar decomposition, the iteration
(\ref{eqn:newton1})--(\ref{eqn:newton2}) can be derived by applying Newton's method
to a matrix equation, here $X^2-A=0$; see~\cite[pp.~139--140]{HigBook08} for details.
As shown in~\cite[Theorem~6.9]{HigBook08}, the iterates $X_k$ converge (under some
natural assumptions) to the principal square root $A^{1/2}$, and they satisfy the
convergence bound
$$\|A^{1/2}-X_{k+1}\|\leq \frac12\,\|X_k^{-1}\|\,\|A^{1/2}-X_k\|^2,$$
where $\|\cdot\|$ denotes any submultiplicative matrix norm.
This is a bound of the classical form (\ref{eqn:classical}) that
shows quadratic convergence.

\subsection{A convergence theorem for the iteration (\ref{eqn:newton1})--(\ref{eqn:newton2})}

The approach for analyzing the iteration (\ref{eqn:newton1})--(\ref{eqn:newton2})
using Theorem~\ref{thm:PotPta} contains some differences in comparison with the one
for the matrix polar decomposition. These differences will be pointed out in the
derivation. Ultimately they will lead to a convergence result giving formally the
same bound as in Theorem~\ref{thm:pp} (with the appropriate modifications), but
containing a strong assumption on the initial matrix $X_0$.

As in Section~\ref{sec:polar} we consider $E={\mathbb C}^{n\times n}$ and
the spectral norm $\|\cdot\|$. Now the mapping $G\,:\,D(G)\rightarrow E$ is defined by
$$G(X)\equiv\frac12 (X+X^{-1}A);\quad\mbox{cf. (\ref{eqn:newton2}).}$$
We need to determine $D(G)$ and define the sets $Z(t)\subseteq D(G)$ for some interval
$T=(0,\tau)$, where we will again specify $\tau=\infty$ in the end.

Clearly, $D(G)$ is equal to the set of the invertible matrices in ${\mathbb C}^{n\times n}$,
so that each $X\in Z(t)$ must be invertible.
The first major difference in comparison to the iteration for the polar decomposition,
where the iterates converge to a unitary matrix,
is that in case of (\ref{eqn:newton1})--(\ref{eqn:newton2}) the conditioning of the
iterates may deteriorate during the process. In fact, when $A$ is (close to) singular, any of
its square roots will be (close to) singular as well. As originally suggested by Pt\'ak
in~\cite{Pta76a,Pta76b} (also cf.~\cite[pp.~20--22]{PotPtaBook84}), such deterioration may be dealt
with by requiring that each $X\in Z(t)$ satisfies
$$\sigma_{\min}(X)\geq h(t),$$
where $h$ is some positive function on $T$ that will be specified below. Further requirements
for each $X\in Z(t)$ are that $X$ commutes with $A$ (cf. (\ref{eqn:newton1})) and that
$\|X-G(X)\|\leq t$, which is the first condition in (ii) of Theorem~\ref{thm:PotPta}.
We thus define for the given matrix $A\in {\mathbb C}^{n\times n}$ and $t\in T=(0,\tau)$
the set
$$Z(t) \equiv\left\{ X\in {\mathbb C}^{n\times n}\,|\, AX=XA,\;\;\sigma_{\min}(X)\geq h(t),\;\;
\|X-G(X)\|\leq t\right\}\;\subset\;D(G).$$
Satisfying (i) and the second condition in (ii) of Theorem~\ref{thm:PotPta} will now prove
the convergence of Newton's method (\ref{eqn:newton1})--(\ref{eqn:newton2}) to a matrix
$X_*\in Z(0)$ with the error norm in each step of the iteration bounded as in (2.c) of
Theorem~\ref{thm:PotPta}. We will have a closer look at the set $Z(0)$ below.

For the second condition in (ii), let $t\in T$ and $X\in Z(t)$ be given.
Since $X$ commutes with $A$, we see that the matrix $\widehat{X}\equiv G(X)$
commutes with $A$. We need to show that
$$\sigma_{\min}(\widehat{X})\geq h(\omega(t))\quad\mbox{and}\quad
\|X-\widehat{X}\|\leq \omega(t)$$
for a positive function $h$ on $T$ and a rate of convergence $\omega$ on $T$. Using
$\|X-\widehat{X}\|\leq t$ and $\sigma_{\min}(X)\geq h(t)$,
we get
\begin{eqnarray*}
\sigma_{\min}(\widehat{X}) &=&\sigma_{\min}(X-(X-\widehat{X}))\geq \sigma_{\min}(X)-\|X-\widehat{X}\|
\geq h(t)-t,
\end{eqnarray*}
from which we obtain the condition
\begin{equation}\label{eqn:cond1}
h(t)-t \geq h(\omega(t)).
\end{equation}
Next, note that $\widehat{X}=G(X)$ can equivalently be written as
$X^2-2 X\widehat{X}+A=0.$
Using this and the fact that $X$ and $\widehat{X}$ commute, we get
\begin{equation}\label{eqn:Asquare}
\widehat{X}^2-A = \widehat{X}^2-A + (X^2-2 X\widehat{X}+A)=(X-\widehat{X})^2.
\end{equation}
Consequently,
\begin{eqnarray}
\|\widehat{X}-G(\widehat{X})\|&=& 
\frac12 \|\widehat{X}^{-1}(A-\widehat{X}^2)\|
\leq \frac12 \|\widehat{X}^{-1}\|\,\|\widehat{X}^2-A\|
= \frac12 \frac{\|(X-\widehat{X})^2\|}{\sigma_{\min}(\widehat{X})}\label{eqn:ineq}\\
&\leq& \frac12 \frac{t^2}{h(t)-t},\nonumber
\end{eqnarray}
which gives the condition
\begin{equation}\label{eqn:cond2}
\frac12 \frac{t^2}{h(t)-t}\leq \omega(t).
\end{equation}
The inequality in (\ref{eqn:ineq}) represents the second major difference in the derivation
compared to the matrix polar decomposition. In (\ref{eqn:ineq}) we have
used the submultiplicativity of the matrix norm, while in (\ref{eqn:eq})--(\ref{eqn:pp3})
no inequalities occur because the matrices are effectively diagonal when considered under the
unitarily invariant norm. The absence of inequalities in the derivation ultimately led to
a convergence bound in Theorem~\ref{thm:pp} that is attained in every step $k\geq 1$.

The two (sufficient) conditions (\ref{eqn:cond1}) and (\ref{eqn:cond2}) for
$\widehat{X}\in Z(\omega(t))$ form a system of functional
inequalities that was ingeniously solved by Pt\'ak in~\cite{Pta76a}. Here it suffices to
say that for any $\gamma\geq 0$ the rate of convergence $\omega$ on $T\equiv (0,\infty)$
and the corresponding function $\sigma$ given by
\begin{equation}\label{eqn:rate}
\omega(t)\equiv \frac12\,\frac{t^2}{(t^2+\gamma^2)^{1/2}}\quad\mbox{and}\quad
\sigma(t)\equiv t-\gamma+(t^2+\gamma^2)^{1/2}
\end{equation}
(i.e. (\ref{eqn:omegaPP}) and (\ref{eqn:sigmaPP}) for $\gamma=1$) and the function
$$h(t)\equiv \gamma+\sigma(t),$$
which is positive on $T=(0,\infty)$, satisfy the conditions (\ref{eqn:cond1}) and
(\ref{eqn:cond2}).

It remains to satisfy the condition (i) in Theorem~\ref{thm:PotPta}, i.e., to show that for some
$t_0>0$ there exists an $X_0\in Z(t_0)$. This will determine the parameter $\gamma\geq 0$, which
so far is not specified. We require that $X_0$ commutes with $A$, and that
$$\sigma_{\min}(X_0)\geq h(t_0)=t_0+(t_0^2+\gamma^2)^{1/2}\quad\mbox{and}\quad
\frac12 \|X_0^{-1}A-X_0\|\leq t_0.$$
If we define
$$t_0\equiv \frac12 \|X_0^{-1}A-X_0\|,$$
then the second inequality is automatically satisfied.
Some simplifications of the first inequality lead to
$$\sigma_{\min}(X_0)\,(\sigma_{\min}(X_0)-2t_0)\geq \gamma^2.$$
Hence we must have $\sigma_{\min}(X_0)-2t_0\geq 0$, and
\begin{equation}\label{eqn:gammadef}
\gamma=\gamma_0\equiv \sigma_{\min}^{1/2}(X_0)(\sigma_{\min}(X_0)-2t_0)^{1/2}
\end{equation}
is a feasible choice of the parameter in order to satisfy all conditions.
Theorem~\ref{thm:PotPta} now shows that under these conditions the
iterative algorithm (\ref{eqn:newton1})--(\ref{eqn:newton2}) is meaningful
and converges to some $X_*\in Z(0)\subseteq \overline{D(G)}$.

Let us now consider the set $Z(0)$. First note that, for any $t>0$ and $X\in Z(t)$,
$$\|\widehat{X}^2-A\|\leq \|X-\widehat{X}\|^2\leq t^2;$$
see (\ref{eqn:Asquare}). Taking the limit $t\rightarrow 0$ we see that each $X\in Z(0)$
must satisfy $X^2=A$. Moreover, each $X\in Z(0)$ commutes with $A$ and
satisfies $\sigma_{\min}(X)\geq h(0)=\gamma_0$. If $\gamma_0>0$, then all
matrices in the limiting set $Z(0)$ are invertible and $Z(0)\subset D(G)$.
If $\gamma_0=0$, then a singular limiting matrix is possible and in this case
$Z(0)$ is not a subset of $D(G)$. This will be discussed further below and numerically
illustrated in the next section.

In summary, we have shown the following result.

\medskip
\begin{theorem}\label{thm:main}
Let $A\in {\mathbb C}^{n\times n}$ and let $X_0\in {\mathbb C}^{n\times n}$
be invertible and commute with $A$. If
\begin{equation}\label{eqn:condition}
\sigma_{\min}(X_0)\geq \|X_0^{-1}A-X_0\|\equiv 2t_0,
\end{equation}
then the iterative algorithm (\ref{eqn:newton1})--(\ref{eqn:newton2})
is meaningful and its sequence of iterates $\{X_k\}$ converges
to a matrix $X_*\in {\mathbb C}^{n\times n}$ that
commutes with $A$ and satisfies $X_*^2=A$. Moreover,
with $\omega$ and $\sigma$ as in (\ref{eqn:rate}), where
$\gamma=\gamma_0$ as in (\ref{eqn:gammadef}), we have
\begin{equation}\label{eqn:bound}
\|X_*-X_k\|\leq \sigma(\omega^{(k)}(t_0)),\quad k=0,1,2,\dots\,.
\end{equation}
If $\gamma_0>0$, then $X_*$ is invertible with $\sigma_{\min}(X_*)\geq \gamma_0$.
\end{theorem}

\medskip
We will now discuss the condition (\ref{eqn:condition}). First note that from
$$\|X_0^{-1}A-X_0\|\leq \|X_0^{-1}\|\,\|A-X_0^2\|=\sigma_{\min}^{-1}(X_0)\,\|A-X_0^2\|$$
it follows that (\ref{eqn:condition}) is implied by
$ \sigma_{\min}^2(X_0)\geq \|A-X_0^2\|$. It therefore appears that the condition is
satisfied whenever $X_0$ commutes with $A$, is well conditioned and sufficiently
close (in norm) to a square root of $A$. However, as we will see below, the condition
(\ref{eqn:condition}) seriously limits the applicability of Theorem~\ref{thm:main}.

The fact that the initial matrix $X_0$ has to satisfy some conditions (unlike in the
case of Theorem~\ref{thm:pp}) is to be expected, since not every matrix $A$ has
a square root; cf.~\cite[Section~1.5]{HigBook08}. If an invertible matrix $X_0$ that
commutes with $A$ and satisfies (\ref{eqn:condition}) can be found, then the convergence
of the Newton iteration guarantees the existence of a square root of $A$.  Pt\'ak referred
to this as an {\em iterative existence proof}~\cite{Pta76b}.

Now let $A\neq 0$ be real and symmetric with the real eigenvalues $\lambda_1\leq\cdots\leq\lambda_n$
and let $X_0=\alpha I$ for some real $\alpha\neq 0$. Then $X_0$ is invertible and
commutes with $A$, and (\ref{eqn:condition}) becomes
\begin{equation}\label{eqn:condalpha1}
\alpha^2\geq \|A-\alpha^2 I\| = \max_{1\leq j\leq n}\,|\lambda_j-\alpha^2| =
\max\,\{\,|\lambda_1-\alpha^2|,\,|\lambda_n-\alpha^2|\,\}.
\end{equation}
This condition cannot be satisfied if $\lambda_1<0$. (Note that in this case the principal 
square root of $A$ does not exist~\cite[Theorem~1.29]{HigBook08}.) 
On the other hand, if $\lambda_1\geq 0$
and hence $A$ is positive semidefinite, then a simple computation shows that each matrix
of the form
\begin{equation}\label{eqn:X0}
X_0=\alpha I\quad\mbox{with}\quad |\alpha|\geq \left(\frac{\lambda_1+\lambda_n}{2}\right)^{1/2}>0
\end{equation}
satisfies (\ref{eqn:condition}) and thus the conditions of Theorem~\ref{thm:main}. If $\lambda_1=0$
and hence $A$ is singular, then for every $|\alpha|\geq (\lambda_n/2)^{1/2}$ equality holds in
(\ref{eqn:condalpha1}). We then get $\gamma_0=0$ in (\ref{eqn:gammadef}), 
and in this case the rate of convergence is
\begin{equation}\label{eqn:singular}
\omega(t)=\frac{t}{2},\quad\mbox{giving}\quad \omega^{(k)}(t)=\frac{t}{2^k}\quad \mbox{and}\quad
\sigma(\omega^{(k)}(t))=\frac{t}{2^{k-1}}.
\end{equation}
A numerical example with a real, symmetric and singular matrix is shown in the right part 
of Figure~\ref{fig:1}.

The following result shows that, despite the inequality in (\ref{eqn:ineq}),
the bound (\ref{eqn:bound}) is attained in every step for certain $A$ and $X_0$.

\medskip
\begin{corollary}\label{cor:exactSQ}
If $A\in {\mathbb R}^{n\times n}$ is symmetric with  eigenvalues
$0\leq \lambda_1\leq\cdots\leq\lambda_n$ and $X_0$ is as in (\ref{eqn:X0}),
then the conditions of Theorem~\ref{thm:main} are satisfied, the sequence
of iterates $\{X_k\}$ of (\ref{eqn:newton1})--(\ref{eqn:newton2}) converges 
to the principal square root $A^{1/2}$ of $A$, and the bound (\ref{eqn:bound}) 
is an equality in each step $k\geq 0$.
\end{corollary}

\medskip
\begin{proof}
We already know that $A$ and $X_0$ satisfy the conditions of Theorem~\ref{thm:main}.
Let $A=Q\Lambda Q^T$ be an eigendecomposition of $A$ with $Q\in {\mathbb R}^{n\times n}$
orthogonal and $\Lambda={\rm diag}(\lambda_1,\dots,\lambda_n)$. Then
$A^{1/2}=Q\Lambda^{1/2} Q^T$ with
$\Lambda^{1/2}={\rm diag}(\lambda_1^{1/2},\dots,\lambda_n^{1/2})$ is the principal square
root of $A$. It suffices to show that $\|A^{1/2}-X_0\|=\sigma(t_0)$, then the result
follows from the last part of Theorem~\ref{thm:PotPta}.

Using $X_0=\alpha I$ and $|\alpha|\geq \left(\frac{\lambda_1+\lambda_n}{2}\right)^{1/2}>0$ we get
\begin{eqnarray*}
t_0 &=& \frac12\|\alpha^{-1}\Lambda -\alpha I\|=\frac{1}{2|\alpha|}\max_{1\leq j\leq n}\,|\lambda_j-\alpha^2|
=\frac{\alpha^2-\lambda_1}{2|\alpha|},\\
\gamma_0 &=& \sigma_{\min}(X_0)^{1/2}\,(\sigma_{\min}(X_0)-2t_0)^{1/2} =
(\alpha^2-2|\alpha| t_0)^{1/2} = \lambda_1^{1/2},
\end{eqnarray*}
and thus
\begin{eqnarray*}
\sigma(t_0) &=& t_0-\gamma_0+(t_0^2+\gamma_0^2)^{1/2} = t_0-\gamma_0+(t_0^2+\alpha^2-2|\alpha| t_0)^{1/2}
= t_0-\gamma_0+(|\alpha|-t_0)\\
&=&|\alpha|-\gamma_0 = |\alpha|-\lambda_1^{1/2} = \max_{1\leq j\leq n}\,|\,|\alpha|-\lambda_j^{1/2}\,| =
\|A^{1/2}-X_0\|,
\end{eqnarray*}
where we have used that $|\alpha|\geq t_0$ and
$|\alpha|\geq \left(\frac{\lambda_1+\lambda_n}{2}\right)^{1/2}>\frac{\lambda_1^{1/2}+\lambda_n^{1/2}}{2}$.
\end{proof}

\medskip
Similar to Theorem~\ref{thm:pp}, the sharpness result in the above corollary could also be derived
using results of Potra and Pt\'ak (e.g.,~\cite[Proposition~5.10]{PotPtaBook84}),
when one interprets the iteration (\ref{eqn:newton1})--(\ref{eqn:newton2}) for a real and
symmetric positive semidefinite matrix $A$ and with $X_0=\alpha I$ in terms of $n$ uncoupled
scalar iterations.

Finally, we note that for a real and symmetric positive definite matrix $A$ with
eigenvalues $0<\lambda_1\leq\cdots\leq\lambda_n$ the condition (\ref{eqn:condition})
with $X_0=\alpha A$ and a real $\alpha\neq 0$ becomes
$$\lambda_1 \geq \max_{1\leq j\leq n} \,|\lambda_j-\alpha^{-2}|=
\max\,\{\,|\lambda_1-\alpha^{-2}|,\,|\lambda_n-\alpha^{-2}|\,\}.$$
An elementary computation shows that this condition cannot be satisfied for any
real $\alpha\neq 0$ when $\lambda_n/\lambda_1>3$. If $1<\lambda_n/\lambda_1\leq 3$,
then $\lambda_n=\mu\lambda_1$ for some $1< \mu\leq 3$, and
$\alpha=(\delta\lambda_1)^{-1/2}$ with $\mu-1\leq \delta\leq 2$ is a feasible
choice. Obviously, the conditions on $A$ and $\alpha$ for the applicability of
Theorem~\ref{thm:main} are in this case much more restrictive than for initial matrices of
the form $X_0=\alpha I$. 

\subsection{Numerical examples}\label{sec:numer}
We will now illustrate Theorem~\ref{thm:main} on some numerical examples.
As originally shown by Higham~\cite{Hig86a}, the Newton iteration for the matrix square root
in the form (\ref{eqn:newton1})--(\ref{eqn:newton2}) is numerically unstable; see
also~\cite[Section~6.4]{HigBook08}. For our numerical experiments we have used the
mathematically equivalent ``IN iteration'' (cf.~\cite[Equation (6.20)]{HigBook08}),
originally derived by Iannazzo~\cite{Ian03}:
\begin{eqnarray}
& & \mbox{$X_0$ is invertible and commutes with $A$, and $E_0=\frac12 (X_0^{-1}A-X_0)$.}\label{eqn:in1}\\
& & X_{k+1} = X_k+E_k\quad\mbox{and}\quad
E_{k+1} = -\frac12 E_kX_{k+1}^{-1}E_k,\quad k=0,1,2,\dots\,.\label{eqn:in2}
\end{eqnarray}
This variant is numerically stable and has a limiting accuracy on the order of the unit
roundoff~\cite[Table~6.2]{HigBook08}. In practice one uses scalings for accelerating 
the convergence~\cite[Section~6.5]{HigBook08}; for brief comments on scalings see 
Section~\ref{sec:concl}.

We use three initial matrices of the form
$$X_0=\alpha_jI\quad\mbox{with}\quad
\alpha_j=2^j \left(\frac{\|A\|}{2}\right)^{1/2},\quad j=1,2,3.$$
This choice of $\alpha_j$ is motivated by the argument leading to (\ref{eqn:X0}),
although this argument only applies to real symmetric matrices.
We point out that for none of our example matrices $A$ an initial matrix
of the form $X_0=\alpha_j A$ satisfies the condition (\ref{eqn:condition}), although
the Newton method always converges when started with them (except for the singular matrix
defined below, where $X_0=\alpha_j A$ cannot be used). This demonstrates the severe restriction
of the condition (\ref{eqn:condition}) for the general applicability of Theorem~\ref{thm:main}.

\begin{figure}
\begin{minipage}{0.49\linewidth}
\includegraphics[width=0.99\linewidth]{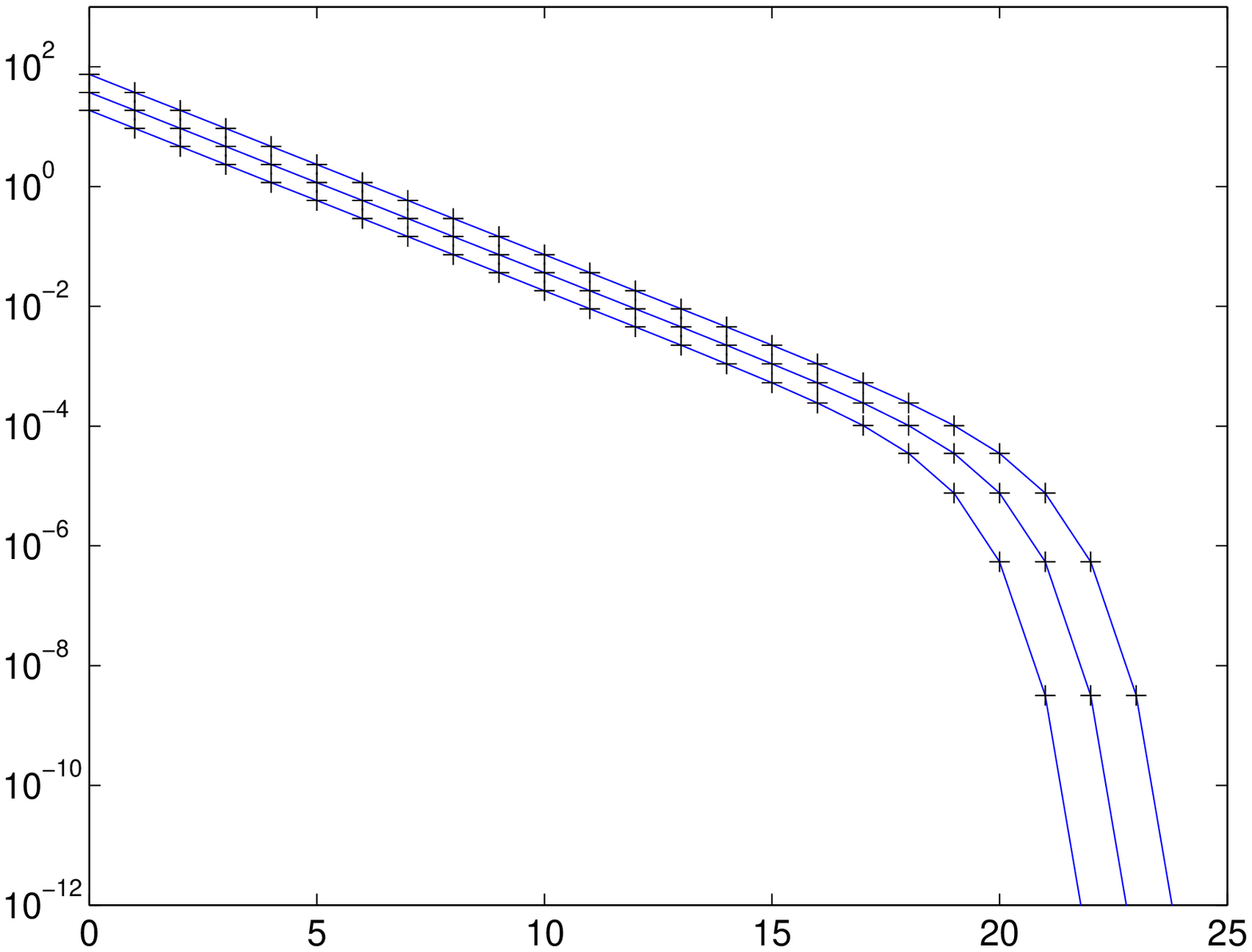}
\end{minipage}
\hfill
\begin{minipage}{0.49\linewidth}
\includegraphics[width=0.99\linewidth]{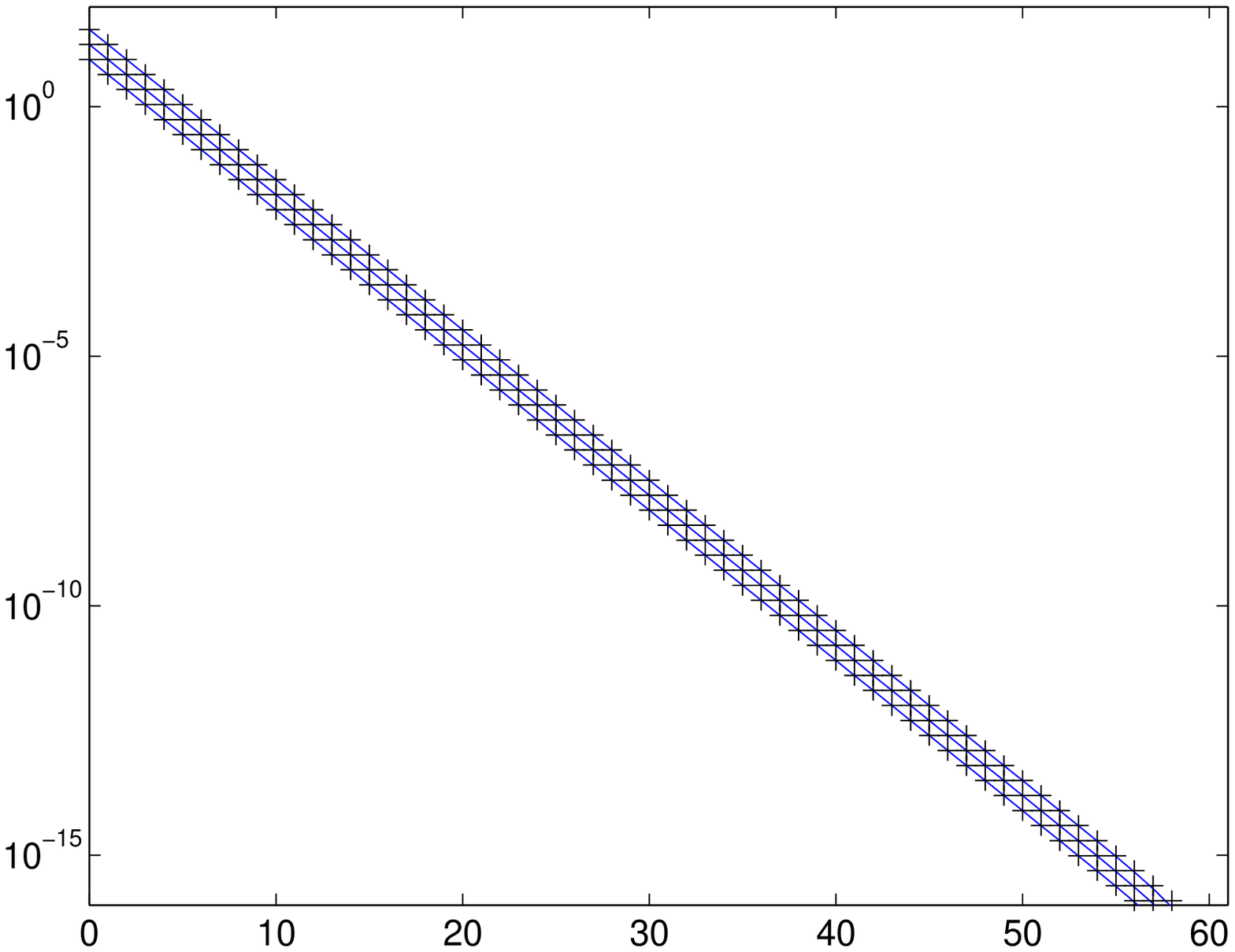}
\end{minipage}
\caption{$\|X_*-X_k\|$ (solid) and $\sigma(\omega^{(k)}(t_0))$ (pluses) for
the Moler matrix (left) and the singular matrix (right).}\label{fig:1}
\end{figure}

\begin{figure}
\begin{minipage}{0.49\linewidth}
\includegraphics[width=0.99\linewidth]{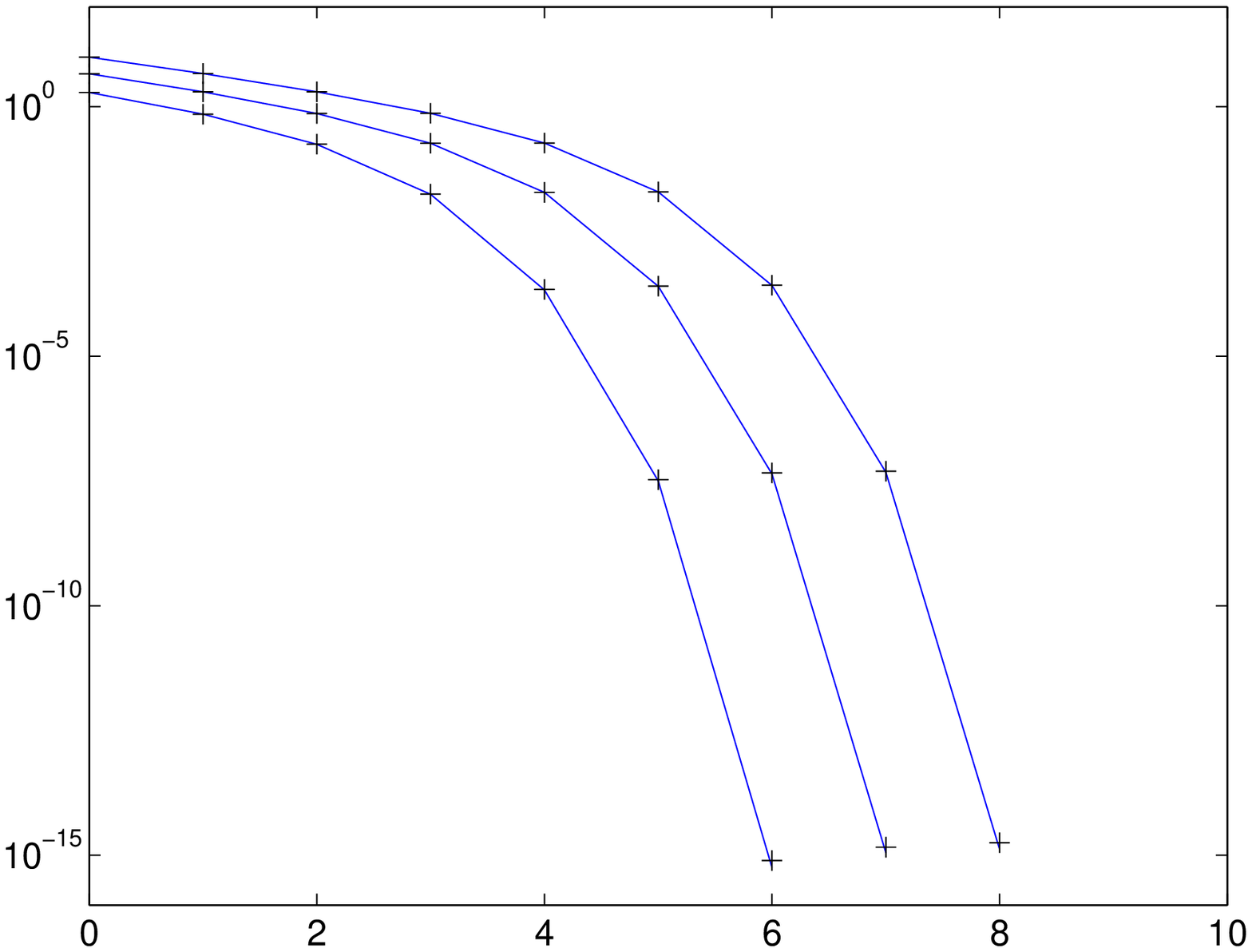}
\end{minipage}
\hfill
\begin{minipage}{0.49\linewidth}
\includegraphics[width=0.99\linewidth]{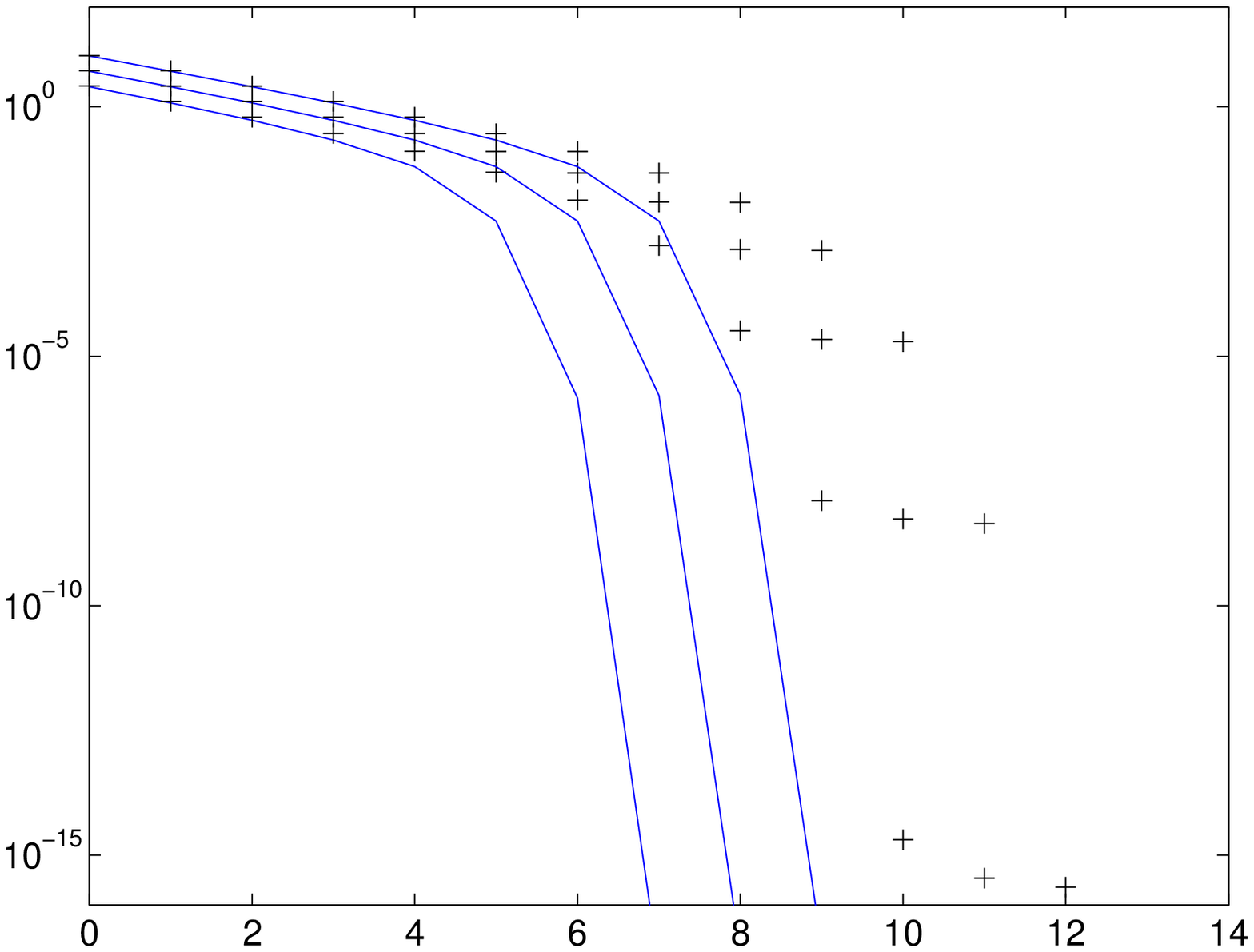}
\end{minipage}
\caption{$\|X_*-X_k\|$ (solid) and $\sigma(\omega^{(k)}(t_0))$ (pluses) for the Jordan blocks
matrix (left) and the modified Jordan blocks matrix (right). }\label{fig:2}
\end{figure}

\begin{figure}
\begin{center}
\includegraphics[width=0.49\linewidth]{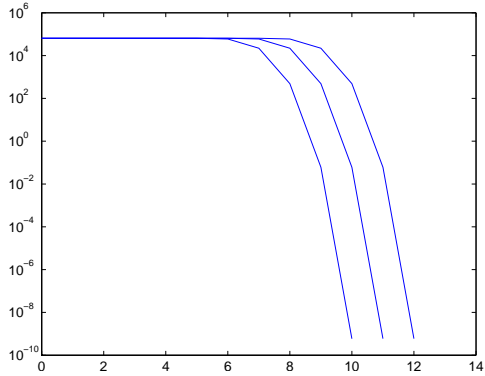}
\caption{$\|X_*-X_k\|$ (solid) for the Frank matrix.}\label{fig:3}
\end{center}
\end{figure}

For the given $X_0$ we first run the iteration (\ref{eqn:in1})--(\ref{eqn:in2}) until it converges
to a matrix $X_*$ that is an approximate square root of $A$. Typically in our examples
$\|A-X_*^2\|\approx 10^{-14}$. We then run the iteration again and compute the error
$\|X_*-X_k\|$ in each step $k\geq 0$. These errors are shown by the solid lines in
Figures~\ref{fig:1}--\ref{fig:3}. The corresponding a priori values $\sigma(\omega^{(k)}(t_0))$
of the upper bound (\ref{eqn:bound}) are shown by the pluses $(+)$ in Figures~\ref{fig:1}--\ref{fig:2}.

We present numerical results for the Moler matrix (left part of Figure~\ref{fig:1}),
the Jordan blocks matrix (left part of Figure~\ref{fig:2}) and the Frank matrix
(Figure~\ref{fig:3}) as defined in Section~\ref{sec:numerPP}. We do not use the 
symmetric indefinite Fiedler matrix from Section~\ref{sec:numerPP} since the Newton 
iteration (\ref{eqn:newton1})--(\ref{eqn:newton2}) for this matrix does not converge 
when started with $X_0=\alpha_j I$. Instead, we present results for
the following matrices:

\smallskip
{\em Singular matrix.\/} The $40\times 40$ diagonal matrix with eigenvalues $0,1,\dots,39$.
The results for this matrix are shown in the right part of Figure~\ref{fig:1}.

\smallskip
{\em Modified Jordan blocks matrix.\/} The $100\times 100$ matrix ${\rm diag}(J_{50}(1.0),J_{50}(2.5))$,
i.e., the previous Jordan blocks matrix with the eigenvalue 1.5 replaced by 1.0. The results
for this matrix are shown in the right part of Figure~\ref{fig:2}.

\smallskip
For the Moler matrix and the singular matrix we have
$\alpha_j \geq \left(\frac{\lambda_1+\lambda_n}{2}\right)^{1/2}$ for $j=1,2,3$.
Hence the bound (\ref{eqn:bound}) is sharp by Corollary~\ref{cor:exactSQ},
which is illustrated by the numerical results. While for the Moler matrix we observe the typical
behavior of Newton's method, the convergence for the singular matrix is linear throughout,
which is explained by (\ref{eqn:singular}). The bound (\ref{eqn:bound}) is very tight
for the Jordan blocks matrix, but for the modified Jordan blocks matrix the actual
convergence acceleration in the quadratic phase is larger than the one predicted by
the bound.

Finally, for the Frank matrix none of the initial matrices $X_0=\alpha_j I$
with $\alpha_j=2^j \left(\|A\|/2\right)^{1/2}$ (here $\|A\|\approx 79.7$)
satisfies the condition (\ref{eqn:condition}), and hence no convergence bound
is given by Theorem~\ref{thm:main}. The constants in condition (\ref{eqn:condition})
have the following values:

\smallskip
\begin{center}
\begin{tabular}{c|r|c}
$j$ & $\sigma_{\min}(X_0)$ & $2t_0=\|X_0^{-1}A-X_0\|$ \\\hline
1 & 9.7710  & 10.1148\\
2 & 19.5420 & 19.7128\\
3 & 39.0839 & 39.1692
\end{tabular}
\end{center}

\smallskip
\noindent Although Theorem~\ref{thm:main} is not applicable, the Newton iteration (eventually)
converges, as shown in Figure~\ref{fig:3}. The figure also shows that the convergence is
characterized by an initial phase of complete stagnation. It is a challenging task
to determine a rate of convergence that can closely describe such behavior.

\section{Discussion and outlook}\label{sec:concl}

The main purpose of Pt\'ak's method of nondiscrete induction is to derive error bounds for iterative
processes that can be tight throughout the iteration. As demonstrated for the Newton
iteration for the matrix polar decomposition and the matrix square root, the approach
can indeed give very tight, even sharp bounds in practical applications. However,
we have also seen that the approach has some limitations, in particular for
(highly) nonnormal matrices. The results for the matrix square root iteration applied to the Frank matrix
illustrate an observation we made: For a general nonnormal matrix it appears to be very
difficult to determine (easily and a priori) an initial matrix $X_0$ that satisfies the assumptions of
Theorem~\ref{thm:main}. The derivation of strategies for computing such $X_0$ is beyond the scope
of this paper. We also did not analyze more closely the observations about the tightness
of the bound (\ref{eqn:bound}) for the two Jordan blocks matrices.
Such investigations are subject of future work.
Note that the difficulty of nonormality does not appear in case of the Newton
iteration for the matrix polar decomposition, since it effectively acts on diagonal matrices.

Convergence of an iteration means that a solution of a certain equation exists. In case of the
Newton iteration for the matrix square root this equation is $X^2-A=0$ for the given (square) matrix $A$.
Finding an initial matrix $X_0$ that satisfies the assumptions of Theorem~\ref{thm:main} guarantees
that the matrix $A$ has a square root (``iterative existence proof''). The fact that existence and
convergence are studied and proven simultaneously may be one reason why the (sufficient) conditions
in Theorem~\ref{thm:main} are more restrictive than other types of conditions for convergence of the
iteration under the assumption that a square root exists. From a practical point of view it would be
of interest to obtain variants of Theorem~\ref{thm:main} with weaker conditions on $X_0$.

As described in Section~\ref{sec:method}, an essential constructive idea in the method
of nondiscrete induction is the comparison of consecutive iterates. The method
therefore is suited best for iterations where a new iterate is generated using information
only from the previous step, and Newton-type methods
are ideal examples in this respect. Methods that generate new iterates based on
information from several or all previous steps possibly can be analyzed using the concept of
{\em multidimensional rates of convergence} described in~\cite[Chapter~3]{PotPtaBook84}.
It also would be interesting to appropriately generalize the theory of Potra and Pt\'ak to
scaled or, more generally, parameter-dependent iterations. As briefly mentioned in
Sections~\ref{sec:numerPP} and~\ref{sec:numer}, scalings are used in practice for
numerical stability properties or for accelerating 
Newton-type matrix iterations. Including parameters
that are not known a priori within the method of nondiscrete induction is an interesting
subject of further work. In the framework established by Theorem~\ref{thm:PotPta} it may
be possible to include parameters depending only on the latest iterate $x_k$ by appropriately
defining the sets $Z(t)$. More sophisticated parameter-dependence requires extensions of
Theorem~\ref{thm:PotPta}. Such extensions may also lead to new convergence results for
nonlinear iterative schemes like Krylov subspace methods (see, e.g.,~\cite{LieStrBook13}),
where new iterates are generated using information from all previous steps and parameters
that are not known a priori.

\medskip
{\bf Acknowledgments.} Many thanks to Nick Higham, G{\'e}rard Meurant,
Zden{\v e}k Strako{\v s} and Petr Tich\'y for helpful comments.

\end{document}